\documentclass[10pt]{amsart}

\usepackage{amsmath,amsthm,amsfonts,amscd,amssymb}
\usepackage{hyperref,wasysym}
\usepackage{verbatim}

\newtheorem{thm}{Theorem}[section]
\newtheorem{cor}[thm]{Corollary}
\newtheorem{lem}[thm]{Lemma}

\theoremstyle{definition}

\theoremstyle{remark}

\begin{document}

\title[Diophantine avoidance and small-height primitive elements]{Diophantine avoidance and small-height primitive elements in ideals of number fields}

\author{Lenny Fukshansky}\thanks{Fukshansky was partially supported by the Simons Foundation grant \#519058}
\author{Sehun Jeong}

\address{Department of Mathematics, 850 Columbia Avenue, Claremont McKenna College, Claremont, CA 91711}
\email{lenny@cmc.edu}
\address{Institute of Mathematical Sciences, Claremont Graduate University, Claremont, CA 91711}
\email{sehun.jeong@cgu.edu}


\begin{abstract} 
Let $K$ be a number field of degree $d$. Then every ideal $I$ in the ring of integers ${\mathcal O}_K$ contains infinitely many primitive elements, i.e. elements of degree $d$. A bound on the smallest height of such an element in $I$ follows from some recent developments in the direction of a 1998 conjecture of W. Ruppert. We prove an explicit bound on the smallest height of such a primitive element in the case of quadratic fields. Further, we consider primitive elements in an ideal outside of a finite union of other ideals and prove a bound on the height of a smallest such element. Our main tool is a result on points of small norm in a lattice outside of an algebraic hypersurface and a finite union of sublattices of finite index, which we prove by blending two previous Diophantine avoidance results. We also obtain a bound for small-norm lattice points in the positive orthant in $\mathbb{R}^d$ with avoidance conditions and use it to obtain a small-height totally positive primitive element in an ideal of a totally real number field outside of a finite union of other ideals. Additionally, we use our avoidance method to prove a bound on the Mahler measure of a generating non-sparse polynomial for a given number field. Finally, we produce a bound on the height of a smallest primitive generator for a principal ideal in a quadratic number field. 
\end{abstract}

\maketitle

\def\A{{\mathcal A}}
\def\B{{\mathcal B}}
\def\C{{\mathcal C}}
\def\D{{\mathcal D}}
\def\F{{\mathcal F}}
\def\x{{\mathcal H}}
\def\I{{\mathcal I}}
\def\J{{\mathcal J}}
\def\K{{\mathcal K}}
\def\L{{\mathcal L}}
\def\M{{\mathcal M}}
\def\N{{\mathcal N}}
\def\O{{\mathcal O}}
\def\R{{\mathcal R}}
\def\s{{\mathcal S}}
\def\V{{\mathcal V}}
\def\W{{\mathcal W}}
\def\X{{\mathcal X}}
\def\Y{{\mathcal Y}}
\def\H{{\mathcal H}}
\def\Z{{\mathcal Z}}
\def\OO{{\mathcal O}}
\def\BB{{\mathbb B}}
\def\cee{{\mathbb C}}
\def\EE{{\mathbb E}}
\def\Nn{{\mathbb N}}
\def\pee{{\mathbb P}}
\def\que{{\mathbb Q}}
\def\real{{\mathbb R}}
\def\zed{{\mathbb Z}}
\def\hyp{{\mathbb H}}
\def\aa{{\mathfrak a}}
\def\HH{{\mathfrak H}}
\def\qbar{{\overline{\mathbb Q}}}
\def\eps{{\varepsilon}}
\def\ahat{{\hat \alpha}}
\def\bhat{{\hat \beta}}
\def\gt{{\tilde \gamma}}
\def\h{{\tfrac12}}
\def\be{{\boldsymbol e}}
\def\bei{{\boldsymbol e_i}}
\def\bff{{\boldsymbol f}}
\def\ba{{\boldsymbol a}}
\def\bb{{\boldsymbol b}}
\def\bc{{\boldsymbol c}}
\def\bm{{\boldsymbol m}}
\def\bk{{\boldsymbol k}}
\def\bi{{\boldsymbol i}}
\def\bl{{\boldsymbol l}}
\def\bq{{\boldsymbol q}}
\def\bu{{\boldsymbol u}}
\def\bt{{\boldsymbol t}}
\def\bs{{\boldsymbol s}}
\def\bv{{\boldsymbol v}}
\def\bw{{\boldsymbol w}}
\def\bx{{\boldsymbol x}}
\def\bX{{\boldsymbol X}}
\def\bz{{\boldsymbol z}}
\def\bwy{{\boldsymbol y}}
\def\bY{{\boldsymbol Y}}
\def\bL{{\boldsymbol L}}
\def\baa{{\boldsymbol\alpha}}
\def\bbb{{\boldsymbol\beta}}
\def\bet{{\boldsymbol\eta}}
\def\bxi{{\boldsymbol\xi}}
\def\bo{{\boldsymbol 0}}
\def\bol{{\boldkey 1}_L}
\def\ep{\varepsilon}
\def\p{\boldsymbol\varphi}
\def\q{\boldsymbol\psi}
\def\rank{\operatorname{rank}}
\def\aut{\operatorname{Aut}}
\def\lcm{\operatorname{lcm}}
\def\sgn{\operatorname{sgn}}
\def\spn{\operatorname{span}}
\def\md{\operatorname{mod}}
\def\Norm{\operatorname{Norm}}
\def\dim{\operatorname{dim}}
\def\det{\operatorname{det}}
\def\Vol{\operatorname{Vol}}
\def\rk{\operatorname{rk}}
\def\Gal{\operatorname{Gal}}
\def\WR{\operatorname{WR}}
\def\WO{\operatorname{WO}}
\def\GL{\operatorname{GL}}
\def\pr{\operatorname{pr}}
\def\Tr{\operatorname{Tr}}

\tableofcontents

\section{Introduction and statement of main results}
\label{intro}

Diophantine avoidance has been studied by several authors in the recent years. This term refers to effective results on existence of points of bounded size (measured by norm or height, depending on the context) in a given algebraic set avoiding some specified subsets. In particular, there is a variety of Siegel's lemma-type results on small-size points in linear spaces and lattices with different avoidance conditions (e.g., \cite{faltings1}, \cite{lf1}, \cite{lf2}, \cite{lf3}, \cite{gaudron1}, \cite{gaudron2}, \cite{henk_thiel}, \cite{faltings2}). The application of avoidance conditions allows to understand how ``well-distributed" points of bounded size are in a given set: if it is possible to find them outside of some prescribed collection of subsets of the set in question, then it suggests that they are evenly distributed, in some appropriate sense. More specifically, it is often interesting to understand how a prescribed avoidance condition affects the upper bound on the size of the ``smallest" element in question. In this paper, we are interested in further investigating small-size points in lattices with avoidance conditions. The main application of our investigation is to small-height generators of number fields satisfying certain natural avoidance conditions.

To start, let $\EE^d$ be a $d$-dimensional Euclidean space, $\Omega \subset \EE^d$ be a lattice of rank $d \geq 2$. Let $\Delta$ be the determinant of $\Omega$, i.e., volume of its fundamental domain, and $\Lambda_1,\dots,\Lambda_s \subset \Omega$ be sublattices of finite indices $D_1,\dots,D_s \geq 2$, $s \geq 1$, so $\det \Lambda_i = D_i \Delta$ for every $1 \leq i \leq s$. Assume that
$$\Omega \not\subseteq \bigcup_{i=1}^s \Lambda_i,$$
and let $\Lambda = \bigcap_{i=1}^s \Lambda_i$. Then $\Lambda$ is a sublattice of $\Omega$ of index no larger than $D := D_1 \cdots D_s$; see~\cite{gruber_lek} for detailed information on lattices and their properties.

For a vector $\bz \in \EE^d$, let
$$|\bz| = \max_{1 \leq i \leq d} |z_i|$$
be its sup-norm, and define 
$$C_d(T) = \left\{ \bz \in \EE^d : |\bz| \leq T \right\}$$
be the cube of side-length $2T$ centered at the origin in $\EE^d$, $T > 0$. For any full-rank lattice $L \subset \EE^d$, define $\lambda_i(L)$, the $i$-th successive minimum of $L$ with respect to $C_d(1)$ to be
\begin{equation}
\label{sm}
\lambda_i(L) = \min \left\{ T > 0 : \dim_{\real} \spn_{\real} (C_d(T) \cap L)  \geq i \right\},
\end{equation}
for each $1 \leq i \leq d$. A theorem of Henk and Thiel \cite[Theorem~1.2]{henk_thiel} guarantees that there exists $\bx \in \Omega \setminus \bigcup_{i=1}^s \Lambda_i$ such that
\begin{equation}
\label{ht}
|\bx| < \frac{2^d D \Delta}{\lambda_1(\Lambda)^{d-1} \Vol_d(C_d(1))} \left( \sum_{i=1}^s \frac{1}{D_i} - \frac{s - 1}{D} + \frac{\lambda_1(\Lambda)^d}{D\Delta} \right),
\end{equation}
where $\Vol_d$ stands for the $d$-dimensional measure on $\EE^d$. This result was obtained using a careful analysis and volume computations on the torus group $\real^d/\Lambda$.

On the other hand, let $P(x_1,\dots,x_d) \in \real[x_1,\dots,x_d]$ be a nonzero polynomial of degree $m$. Let $S_1,\dots,S_d$ be finite subsets of $\zed$ with $|S_i| \geq m+1$ for each $1 \leq i \leq d$ and let $\bv_1,\dots,\bv_d$ be linearly independent vectors in our lattice $\Omega$. Then Theorem~4.2 of~\cite{lf3} implies that there exist coefficients $\xi_i \in S_i$ for $1 \leq i \leq d$ so that
\begin{equation}
\label{p_avoid}
P \left( \sum_{i=1}^d \xi_i \bv_i \right) \neq 0.
\end{equation}
The proof of this theorem (Theorem~4.2 of~\cite{lf3}) is an application of Alon's Combinatorial Nullstellensatz~\cite{alon}. Our first goal is to bridge the two above-mentioned avoidance results and prove the following theorem.

\begin{thm} \label{dioph_avoid} Let the notation be as above. Then there exists
$$\bz \in \Omega \setminus \left( \bigcup_{i=1}^s \Lambda_i \right),$$
such that $P(\bz) \neq 0$ and
\begin{equation}
\label{z-bound}
|\bz| \leq \frac{d(D(m+2)+2) D\Delta}{2\lambda_1(\Lambda)^{d-1}} \max \left\{ 1, \frac{2^d}{\Vol_d(C_d(1))} \left( \sum_{i=1}^s \frac{1}{D_i} - \frac{s - 1}{D} + \frac{\lambda_1(\Lambda)^d}{D\Delta} \right) \right\}.
\end{equation} 
\end{thm}

\noindent
We prove Theorem~\ref{dioph_avoid} in Section~\ref{DA} using~\eqref{ht}, \eqref{p_avoid} and Minkowski's Successive Minima Theorem. 
\smallskip

We apply our Theorem~\ref{dioph_avoid} in the context of algebraic number fields. Let $K$ be a number field of degree $d = [K:\que] \geq 1$. Let $\sigma_1,\dots,\sigma_d : K \hookrightarrow \cee$ be the embeddings of $K$, ordered so that the first $r_1$ of them are real and the remaining $2r_2$ are conjugate pairs of complex embeddings so that $\sigma_{r_2+j} = \bar{\sigma}_j$ for 
$$d = r_1 + 2r_2.$$

\noindent
An element $\alpha \in K$ is called {\it primitive} if $K = \que(\alpha)$. This is equivalent to the condition that $\deg_{\que}(\alpha) = d$, and hence there are infinitely many primitive elements in $K$. A conjecture of Ruppert \cite{ruppert} asserts that there exists a primitive element $\alpha \in K$ such that
$$h(\alpha) \leq c(d) |\Delta_K|^{\frac{1}{2d}},$$
where $h$ is the absolute Weil height, $\Delta_K$ is the discriminant of the number field $K$, and $c(d)$ is a constant depending only on the degree $d$; we review the height machinery in Section~\ref{heights}. Ruppert himself proved this conjecture for quadratic number fields and for totally real fields of prime degree. There has been quite a bit of later work on this conjecture; for instance, Vaaler and Widmer \cite{vaaler_widmer-1} proved the conjecture for number fields with at least one real embedding. More generally, a slightly weaker bound is obtained by Pazuki and Widmer in \cite[Lemma~7.1]{widmer}:
\begin{equation}
\label{wd0}
h(\alpha) \leq |\Delta_K|^{\frac{1}{d}}.
\end{equation}
In fact, a more detailed result follows from Lemma~7.1 of~\cite{widmer}. Let $\O_K \subset K$ be the ring of integers of $K$ and let $I \subseteq \O_K$ be an ideal in this ring. It is not difficult to see that $I$ contains infinitely many primitive elements, and a straight-forward modification of the proof of Lemma~7.1 of~\cite{widmer} (replacing $\O_K$ by $I$) produces a primitive element $\alpha \in I$ with
\begin{equation}
\label{wd1}
h(\alpha) \leq \Nn_K(I)^{\frac{2}{d}} |\Delta_K|^{\frac{1}{d}},
\end{equation}
where $\Nn_K(I) = |\O_K/I|$ is the norm of the ideal $I$. To this end, we obtain a more concrete result in the case of quadratic number fields.

\begin{thm} \label{quad_bound} Let $D$ be a squarefree integer and let $K=\que(\sqrt{D})$ be a quadratic number field. Let $I \subseteq \O_K$ be an ideal with the canonical basis $\{ a, b+g\delta \}$, as described in~\eqref{abg}, so that $b < a$. Let
$$h_{\min}(I) = \min \left\{ h(\alpha) : \alpha \in I, K = \que(\alpha) \right\}.$$
If $D \not\equiv 1 (\md 4)$, then
$$\sqrt{ag} < h_{\min}(I) \leq g \left( \frac{2b+ \sqrt{|\Delta_K|}}{2} \right),$$
and additionally $h_{\min}(I) > g \sqrt{|\Delta_K|}/2$ if $D < 0$.
\smallskip

\noindent
If $D \equiv 1 (\md 4)$, then
$$\sqrt{ag} < h_{\min}(I) \leq g \left( \frac{(2b+1) + \sqrt{|\Delta_K|}}{2} \right),$$
and additionally $h_{\min}(I) > g \sqrt{|\Delta_K|}/2$ if $D < 0$.
\end{thm}

\noindent
To compare the bounds of Theorem~\ref{quad_bound} to that of~\eqref{wd1}, notice that 
$$bg < ag = \left\{ \begin{array}{ll}
\Nn_K(I) & \mbox{if $D \not\equiv 1 (\md 4)$,} \\
2\Nn_K(I) & \mbox{if $D \equiv 1 (\md 4)$.}
\end{array}
\right.$$
For example, in the case $D \not\equiv 1 (\md 4)$, if $D > 4$ and $b > 2$, we obtain
$$h_{\min}(I) \leq g \left( \frac{2b+ \sqrt{|\Delta_K|}}{2} \right) < \frac{bg \sqrt{|\Delta_K|}}{2} < \frac{\Nn_K(I) \sqrt{|\Delta_K|}}{2},$$
which is a slight improvement on~\eqref{wd1}, since $d=2$; there is a similar comparison in the case $D \equiv 1 (\md 4)$.
We review all the necessary notation and prove Theorem~\ref{quad_bound} in Section~\ref{quad}.
\smallskip

To generalize the above setting, our next result produces a small-height primitive element in an ideal of a number field outside of a finite collection of other ideals.

\begin{thm} \label{ideal_miss} Let $K$ be a number field of degree $d \geq 2$ and let $I \subset \O_K$ be an ideal and $J_1,\dots,J_s \subsetneq I$ be distinct ideals in $\O_K$. Let $J = J_1 \cdots J_s$. Then there exists an element $\alpha \in I \setminus \bigcup_{i=1}^s J_i$ such that $K = \que(\alpha)$ and
\begin{eqnarray}
\label{main}
h(\alpha) & \leq & \frac{d\ \Nn_K(J) |\Delta_K|^{1/2} \left( (d^2-d+4)\ \Nn_K(J)  + 4 \Nn_K(I) \right)}{4 \Nn_K(I)} \times \nonumber \\
& & \times \max \left\{ 1, \left( \frac{2}{\pi} \right)^{r_2} \left( \sum_{i=1}^s \frac{\Nn_K(I)}{\Nn_K(J_i)} - \frac{(s-1) \Nn_K(I)}{\Nn_K(J)}  + \frac{\Nn_K(J)^{d-1}}{|\Delta_K|^{1/2}} \right) \right\}.
\end{eqnarray}
\end{thm}

\noindent
We want to emphasize that imposing an additional avoidance condition on an existence result for a primitive element is natural in the context of the Primitive Element Theorem, which asserts the existence of infinitely many primitive elements in a number field. The standard proof of this celebrated theorem (see, for instance, Theorem~2.2 of~\cite{stewart-tall}) produces a primitive element for a number field $K$ by constructing an element in $K$ that avoids a finite collection of linear equations over $\que$. In other words, the avoidance idea is already intrinsic to this construction. In fact, our avoidance argument produces, among other things, an alternate proof of the (effective version of the) Primitive Element Theorem.

Notice that the bound of~\eqref{main} is explicit in terms of the invariants of $K$ and norms of the ideals involved. To get a better grasp on its order of magnitude, it may be helpful to rewrite in a less explicit form as
\begin{equation}
\label{order_bnd}
h(\alpha) \ll_{K,s} \frac{\Nn_K(J)^{d+1}}{\Nn_K(I)}.
\end{equation}
We prove Theorem~\ref{ideal_miss} in Section~\ref{pf_ideal_miss}, where our main tool is Theorem~\ref{dioph_avoid}. In Section~\ref{cone}, we prove an analogous result for primitive totally positive elements in an ideal of a totally real number field outside of a finite union of other ideals (Corollary~\ref{ideal_pos}). This result follows as a consequence of a theorem we obtain on existence of a small-norm positive point in a lattice avoiding a finite union of sublattices of the same rank and an algebraic set (Theorem~\ref{pos_cone}). This theorem is proved by blending together the method of proof of Theorem~\ref{dioph_avoid} with the results on small-height lattice points in positive cones obtained in~\cite{lf_siki}. We also use our avoidance method in Section~\ref{mahler} to prove the existence of a monic polynomial $f(x) \in \zed[x]$ with bounded Mahler measure and all nonzero coefficients so that a given number field $K$ is isomorphic to $\que[x]/\left< f(x) \right>$.
\smallskip

So far, we were concerned with small-height primitive elements, i.e. generators of number fields, contained in a prescribed ideal. Now we turn to small-height generators of ideals in number fields. It is well known that every ideal $I \subseteq \O_K$ can be generated by two elements (see, e.g., Theorem~5.20 of~\cite{stewart-tall}), and principal ideals (i.e., ideals generated by a single element) correspond to the identity element in the class group of $K$. If $\mu \in \O_K$ is a generator for a principal ideal $I \subseteq \O_K$ and $\eps \in \O_K^{\times}$ is a unit, then $\eps \mu$ is another generator for $I$. Dirichlet units theorem states that the rank of the unit group $\O_K^{\times}$ is $r = r_1 + r_2 - 1$ (see, e.g., \cite{stewart-tall} for further details). In particular, $\O_K^{\times}$ is infinite unless $K$ is an imaginary quadratic field. Hence, the same principal ideal $I$ can have infinitely many generators, and while the norm of all of them is the same ($= \Nn_K(I)$), their heights can be very different and arbitrarily large. Theorem~1.1 of~\cite{akhtari_vaaler-1} by Akhtari and Vaaler then implies that there exists a generator $\mu$ for $I$ such that
\begin{equation}
\label{av1}
h(\mu) \leq \Nn_K(I)^{1/d} \prod_{j=1}^r h(\eps_j)^{1/2},
\end{equation}
where $\eps_1,\dots,\eps_r$ are multiplicatively independent units in $\O_K^{\times}$. In order to obtain a bound on $h(\mu)$ in terms of the invariants of $I$ and $K$ alone, we now need a result on existence of small-height multiplicatively independent units in $\O_K^{\times}$. Such a result has been obtained recently also by Akhtari and Vaaler in~\cite{akhtari_vaaler-2}, but with a bound on the product of logarithmic heights of the units in question. Hence, using this bound in conjunction with~\eqref{av1} would essentially require reversing arithmetic-geometric mean inequality. On the other hand, we can obtain the following bound on $h(\mu)$ in the case of a quadratic number field using a different method.

\begin{thm} \label{ideal_gen} Let $D$ be a squarefree integer and let $K=\que(\sqrt{D})$ be a quadratic number field. Let $I \subseteq \O_K$ be a principal ideal with the canonical basis $\{ a, b+g\delta \}$, as described in~\eqref{abg}. Define
\begin{equation}
\label{HI}
H(I) := \left\{ \begin{array}{ll}
\max \left\{  \frac{|a|}{g},\frac{|2b|}{g},\frac{|b^2-D|}{ag}  \right\} & \mbox{if $D \not\equiv 1\ (\md 4)$} \\
\max \left\{ \frac{|2a|}{g}, \frac{2 |2b+g|}{g} ,\frac{|4b^2+4b-Dg+g|}{2ag} \right\} & \mbox{if $D \equiv 1\ (\md 4)$},
\end{array}
\right.
\end{equation}
Then there exists a primitive element $\mu \in I$ such that $I = \left< \mu \right>$ and
$$
h(\mu) \leq \left\{ \begin{array}{ll}
\left( a+b+g\sqrt{|D|} \right) (14 H(I))^{5H(I)} & \mbox{if $D \not\equiv 1\ (\md 4)$} \\
\left( a + \frac{2b+g+ g\sqrt{|D|}}{2} \right) (14 H(I))^{5H(I)} & \mbox{if $D \equiv 1\ (\md 4)$},
\end{array}
\right.
$$
\end{thm}

\noindent
Theorem~\ref{ideal_gen} is essentially the quadratic case of a general difficult problem: given an infinite set of elements of fixed norm $N$ in $\O_K$, prove the existence of an element of bounded height in this set with the bound depending only on $N$ and invariants of the number field $K$. We prove Theorem~\ref{ideal_gen} in Section~\ref{small_gen}. Our main tool is a result of Kornhauser~\cite{kornhauser} on small-height zeros of integral binary quadratic equations with integer coefficients. We are now ready to proceed.

\bigskip

\section{Avoiding sublattices and algebraic sets}
\label{DA}

The goal of this section is to prove Theorem~\ref{dioph_avoid}, establishing an effective result on existence of a point in a lattice avoiding a union of full-rank sublattices and an algebraic hypersurface. Let the notation be as in the statement of the theorem. Then Minkowski's successive minima theorem (see, for instance, \cite[Section~9.1, Theorem 1]{gruber_lek}) implies that $\prod_{i=1}^d \lambda_i(\Lambda) \leq D \Delta$, and so
\begin{equation}
\label{mink}
\lambda_d(\Lambda) \leq \frac{D \Delta}{\lambda_1(\Lambda)^{d-1}},
\end{equation}
since $0 < \lambda_1(\Lambda) \leq \dots \leq \lambda_{d-1}(\Lambda) \leq \lambda_d(\Lambda)$.

Let $\bv_1,\dots,\bv_d \in \Lambda$ be linearly independent vectors corresponding to these successive minima, so $|\bv_i| = \lambda_i(\Lambda)$. Let $\bx \in \Omega \setminus \bigcup_{i=1}^s \Lambda_i$ be a vector satisfying~\eqref{ht}, then there is a subset of $d-1$ vectors among $\bv_1,\dots,\bv_d$ with which $\bx$ is linearly independent, assume these are $\bv_2,\dots,\bv_d$, since they have larger norm. Notice that for any integer $k$ and any vector $\bwy \in \Omega$, $Dk\bwy \in \Lambda_i$ for every $1 \leq i \leq s$, and hence the vector $(Dk+1)\bx \not\in \Lambda_i$ for every $1 \leq i \leq s$. Therefore, for all integers $k,n_2,\dots,n_d$, we have a collection of vectors
\begin{equation}
\label{all_vectors}
(Dk+1)\bx + n_2 \bv_2 + \dots + n_d \bv_d \in \Omega \setminus \bigcup_{i=1}^s \Lambda_i.
\end{equation}
Let $[\ ]$ denote the integer part function, and define the sets
\begin{align}
\label{S1S2}
& S_1 = \{ Dk + 1 : k \in \zed, -[m/2]-1 \leq k \leq [m/2]+1 \}, \nonumber \\
& S_2 = \{ j \in \zed : -[m/2]-1 \leq j \leq [m/2]+1 \},
\end{align}
consisting of $m+1$ integers each. Then~\eqref{p_avoid} implies that there exists a vector 
$$\bxi = (\xi_1,\xi_2,\dots,\xi_d) \in S_1 \times S_2 \times \dots \times S_2$$
so that 
$$P\left( \xi_1 \bx + \sum_{i=2}^d \xi_i \bv_i \right) \neq 0.$$
Let $\bz = \xi_1 \bx + \sum_{i=2}^d \xi_i \bv_i $ for this choice of $\bxi$. We now need to estimate its sup-norm. Since
$$\lambda_1(\Lambda) = |\bv_1| \leq \lambda_2(\Lambda) = |\bv_2| \leq \dots \leq \lambda_d(\Lambda) = |\bv_d|,$$
we can put together~\eqref{ht} and~\eqref{mink} to obtain
\begin{eqnarray*}
|\bz| & \leq & d \left( \max_{1 \leq i \leq d} |\xi_i| \right) \left( \max \{ |\bx|, \lambda_d(\Lambda) \} \right) \\
& \leq & \frac{d(D(m+2)+2) D\Delta}{2\lambda_1(\Lambda)^{d-1}} \max \left\{ 1, \frac{2^d}{\Vol_d(C_d(1))}  \left( \sum_{i=1}^s \frac{1}{D_i} - \frac{s - 1}{D} + \frac{\lambda_1(\Lambda)^d}{D\Delta} \right) \right\}.
\end{eqnarray*}
This completes the proof.

\bigskip

\section{Heights and additional notation}
\label{heights}

In this section, we set up the notation needed for the proofs of our Theorems~\ref{quad_bound}, \ref{ideal_miss} and \ref{ideal_gen}. Notice that the space $K_{\real} := K \otimes_{\que} \real$ can be viewed as a subspace of
$$\left\{ (\bx,\bwy) \in \real^{r_1} \times \cee^{2r_2} : y_{r_2+j} = \bar{y}_j\ \forall\ 1 \leq j \leq r_2 \right\} \cong \real^{r_1} \times \cee^{r_2} \subset \cee^d.$$
Here, in the containment 
$$\real^{r_1} \times \cee^{r_2} = \real \times \dots \times \real \times \cee \times \dots \times \cee \subset \cee \times \dots \times \cee \times \cee \times \dots \times \cee = \cee^d$$
each copy of $\real$ is identified with the real part of the corresponding copy of $\cee$ in which it is contained. Then $K_{\real}$ is a $d$-dimensional Euclidean space with the bilinear form induced by the trace form on $K$:
$$\left< \alpha,\beta \right> := \Tr_K(\alpha \bar{\beta}) \in \real,$$
for every $\alpha,\beta \in K$, where $\Tr_K$ is the number field trace on $K$. We also define the sup-norm on $K_{\real}$ by
$$|\bx| := \max \{ |x_1|,\dots,|x_d| \},$$
for any $\bx \in K_{\real}$, where $|x_j|$ stands for the usual absolute value of $x_j$ on $\cee$. Let $\Sigma_K = (\sigma_1,\dots,\sigma_d) : K \hookrightarrow K_{\real}$ be the Minkowski embedding, then for any ideal $I \subseteq \O_K$ the image $\Sigma_K(I)$ is a lattice of full rank in $K_{\real}$. We define the determinant of a full-rank lattice to be the absolute value of the determinant of any basis matrix for the lattice, then
\begin{equation}
\label{ideal_det}
\det(\Sigma_K(I)) = \Nn_K(I) |\Delta_K|^{1/2},
\end{equation}
as follows, for instance, from Corollary 2.4 of \cite{bayer}. 
\smallskip

Next we normalize absolute values and introduce the standard height function. Let us write $M(K)$ for the set of places of $K$. For each $v \in M(K)$ let $d_v = [K_v : \que_v]$ be the local degree, then for each $u \in M(\que)$, $\sum_{v \mid u} d_v = d$. We select the absolute values so that $|\ |_v$ extends the usual archimedean absolute value on $\que$ when $v \mid \infty$, or the usual $p$-adic absolute value on $\que$ when $v \nmid \infty$. With this choice, the product formula reads
$$\prod_{v \in M(K)} |\alpha|_v^{d_v} = 1,$$
for each nonzero $\alpha \in K$. We define the multiplicative Weil height on algebraic vectors $\baa = (\alpha_1,\dots,\alpha_n) \in K^n$ as
$$h(\baa) = \prod_{v \in M(K)} \max \{1,|\alpha_1|_v,\dots,|\alpha_n|_v \}^{\frac{d_v}{d}},$$
for all $n \geq 1$. This height is absolute, meaning that it is the same when computed over any number field $K$ containing $\alpha_1,\dots,\alpha_n$: this is due to the normalizing exponent~$1/d$ in the definition. Hence we can compute height for points defined over $\qbar$.
\smallskip

We also review a couple useful well-known properties of heights. The first can be found, for instance, as Lemma~2.1 of~\cite{lf3}.

\begin{lem} \label{ht_sum} Let $\xi_1,\dots,\xi_m \in \qbar$ and $\bx,\dots,\bx_m \in \qbar^n$ for $m,n \geq 1$. Then
$$h \left( \sum_{j=1}^m \xi_j \bx_j \right) \leq m h(\bxi) \prod_{j=1}^m h(\bx_j),$$
where $\bxi = (\xi_1,\dots,\xi_m)$.
\end{lem}

Next is Lemma 4.1 of~\cite{lf_siki}: while in that paper the lemma is stated for totally real fields, its proof holds verbatim for any number field with our definition of Minkowski embedding~$\Sigma_K$.

\begin{lem} \label{ht_sigma} For any $\alpha \in \O_K$,
$$1 \leq h(\alpha) \leq |\Sigma_K(\alpha)|,$$
where $|\ |$ stands for the sup-norm on $K_{\real}$, as above.
\end{lem}

\bigskip

\section{Quadratic fields}
\label{quad}

In this section we review the additional necessary notation specific to Theorem~\ref{quad_bound} and prove this theorem. First notice that for any number field $K$ and $\alpha \in \O_K$, 
$$h(\alpha) = \prod_{v \mid \infty} \max \{1, |\alpha|_v\}^{\frac{d_v}{d}} \geq \left( \prod_{v \mid \infty} |\alpha|_v^{d_v} \right)^{\frac{1}{d}} = \left( \prod_{j=1}^d |\sigma_j(\alpha)| \right)^{\frac{1}{d}} = \Nn_K(\alpha)^{\frac{1}{d}}.$$
Now let $D$ be a squarefree integer and $K=\que(\sqrt{D})$ be a quadratic field. Let $I \subseteq \O_K$ be an ideal. Then there exists a unique integral basis $a,b+g\delta$ for $I$, called the canonical basis, where
\begin{equation}
\label{delta}
\delta = \left\{ \begin{array}{ll}
- \sqrt{D} & \mbox{if $K=\que(\sqrt{D})$, $D \not\equiv 1 (\md 4)$} \\
\frac{1-\sqrt{D}}{2} & \mbox{if $K=\que(\sqrt{D})$, $D \equiv 1 (\md 4)$},
\end{array}
\right.
\end{equation}
and $a,b,g \in \zed_{\geq 0}$ such that
\begin{equation}
\label{abg}
b < a,\ g \mid a,b,\text{ and } ag \mid \Nn_K(b+g\delta),
\end{equation}
see Section 6.3 of \cite{buell} for further details. The embeddings $\sigma_1, \sigma_2 : K \to \cee$ are given by
$$\sigma_1(x+y\sqrt{D}) = x+y\sqrt{D},\ \sigma_2(x+y\sqrt{D}) = x-y\sqrt{D}$$
for each $x+y\sqrt{D} \in K$, where $D$ is positive for real quadratic field $K$ and negative for imaginary quadratic field $K$. The number field norm on $K$ is given by
$$\Nn_K(x+y\sqrt{D}) = \sigma_1(x+y\sqrt{D}) \sigma_2(x+y\sqrt{D}) = \left( x+y\sqrt{D} \right) \left( x-y\sqrt{D} \right).$$
The discriminant of $K$ is 
\begin{equation}
\label{disc}
\Delta_K = \left\{ \begin{array}{ll}
4D & \mbox{if $K=\que(\sqrt{D})$, $D \not\equiv 1 (\md 4)$} \\
D & \mbox{if $K=\que(\sqrt{D})$, $D \equiv 1 (\md 4)$,}
\end{array}
\right.
\end{equation}
and the norm of the ideal $I$ with the canonical basis as in~\eqref{abg} above is 
\begin{equation}
\label{ideal_norm}
\Nn_K(I) = \left\{ \begin{array}{ll}
ag & \mbox{if $D \not\equiv 1 (\md 4)$,} \\
ag/2 & \mbox{if $D \equiv 1 (\md 4)$.}
\end{array}
\right.
\end{equation}
Observe that an ideal $I$ with the canonical basis as in~\eqref{abg} above can be written as $I = gJ$ for the corresponding ideal $J = \frac{1}{g} I \subseteq \O_K$, since $g \mid a,b$.  Hence we start by restricting our consideration to ideals with $g=1$. Then the bound in~\eqref{wd1} in the case of a quadratic field can be written as
\begin{equation}
\label{N_h}
\sqrt{\Nn_K(\alpha)} \leq h(\alpha) \ll a \sqrt{|D|}.
\end{equation}
We will show that in this case the power on $\sqrt{D}$ cannot in general be reduced. First observe that an element $\alpha \in I$ is primitive if and only if it is of the form 
$$\alpha = xa+y(b+\delta)$$
with $x,y \in \zed$ and $y \neq 0$.
\medskip

{\it Case 1:} Suppose $D \in \zed$ is squarefree, $D \not\equiv 1\ (\md 4)$ and $K=\que(\sqrt{D})$. Then $K$ is a real quadratic if $D > 0$ and $K$ is an imaginary quadratic if $D<0$. Take an ideal
$$I = \spn_{\zed} \{ a, b-\sqrt{D} \} \subseteq \O_K$$
with $a \mid \Nn_K(b-\sqrt{D}) = |b^2-D|$.  Then
\begin{eqnarray}
\label{N1}
\Nn_K(\alpha) & = & \left| \left( xa + y(b-\sqrt{D}) \right) \left( xa + y(b+\sqrt{D}) \right) \right| \nonumber \\
& = & \left| x^2a^2 + 2xyab + y^2(b^2-D) \right| = a \left| x^2a + 2xyb + y^2 \left( \frac{b^2-D}{a} \right) \right| > a,
\end{eqnarray}
where $(b^2-D)/a \in \zed$. Further,
\begin{equation}
\label{N2}
\Nn_K(\alpha) = \left| x^2a^2 + 2xyab + y^2(b^2-D) \right| = \left| (xa+yb)^2 - y^2D \right| > |D|,
\end{equation}
if $D < 0$. On the other hand,
\begin{align}
\label{h_a_mu}
h(\alpha) & = \prod_{v \in M(K)}\max\{1,|\alpha|_v\}^{\frac{d_v}{2}} = \prod_{v | \infty}\max\{1,|\alpha|_v\}^{\frac{d_v}{2}} = \left( \prod_{j=1}^2 \max \{1,|\sigma_i(\alpha)|\} \right)^{\frac{1}{2}} \nonumber \\
& \leq \frac{1}{2}\left(\max \{1,|\sigma_1(\alpha)|\} + \max \{1,|\sigma_2(\alpha)|\right) \leq \frac{1}{2}\left(2 \cdot \max \{|\sigma_1(\alpha)|,|\sigma_2(\alpha)|\}\right) \nonumber \\
& = \max \{|\sigma_1(\alpha)|,|\sigma_2(\alpha)|\} = \max \left\{\left|(xa+yb)-y\sqrt{D}\right|,\left|(xa+yb)+y\sqrt{D}\right|\right\} \nonumber \\
&\leq |x|a+|y|b+|y|\sqrt{|D|}.
\end{align}
Taking the minimum over all primitive elements $\alpha \in I$, we see that
\begin{align}
\label{N3}
\min \{ h(\alpha) : \alpha \in I, K = \que(\alpha) \} & \leq \min \{ |x|a+|y|b+|y|\sqrt{|D|} : x,y \in \zed, y \neq 0\} \nonumber \\
& \leq b + \sqrt{|D|},
\end{align}
where the last inequality is obtained by taking $x=0$, $y=1$. Putting together \eqref{N_h}, \eqref{N1}, \eqref{N2} and \eqref{N3}, we obtain the $D \not\equiv 1 (\md 4)$ case of Theorem~\ref{quad_bound} in case $g=1$.
\medskip

{\it Case 2:} Suppose $D \in \zed$ is squarefree, $D \equiv 1\ (\md 4)$ and $K=\que(\sqrt{D})$. Again, $K$ is a real quadratic if $D > 0$ and $K$ is an imaginary quadratic if $D<0$. Take an ideal
$$I = \spn_{\zed} \left\{ a, b + \frac{1-\sqrt{D}}{2} \right\} \subseteq \O_K$$
with $a \mid \Nn_K\left( \frac{(2b + 1) - \sqrt{D}}{2} \right) = \frac{|(2b+1)^2 - D|}{4} = \left| b^2 + b - \frac{D-1}{4} \right|$.  Then
\begin{eqnarray}
\label{N1-2}
\Nn_K(\alpha) & = & \left| \left( xa + y \left(b + \frac{1-\sqrt{D}}{2} \right) \right) \left( xa + y \left(b + \frac{1+\sqrt{D}}{2} \right) \right) \right| \nonumber \\
& = & \left| x^2a^2 + (2b+1)axy + y^2 \left( b^2 + b - \frac{D-1}{4} \right) \right| \nonumber \\
& = & a \left| x^2a + (2b+1)xy + \frac{y^2}{a} \left( b^2 + b - \frac{D-1}{4} \right) \right| > a,
\end{eqnarray}
where $\frac{1}{a} \left( b^2 + b - \frac{D-1}{4} \right) \in \zed$. Further,
\begin{equation}
\label{N2-2}
\Nn_K(\alpha) = \left| x^2a^2 + (2b+1)axy + y^2 (b^2 + b +1/4) - y^2D/4 \right| > |D|/4,
\end{equation}
if $D < 0$. On the other hand,
\begin{align}
\label{h_a_mu-1}
h(\alpha) & \leq \max \{|\sigma_1(\alpha)|,|\sigma_2(\alpha)|\} \nonumber \\
& = \max \left\{ \left| xa + y \left(b + \frac{1-\sqrt{D}}{2} \right) \right|, \left| xa + y \left(b + \frac{1+\sqrt{D}}{2} \right) \right| \right\}  \nonumber \\
& \leq |x|a+\frac{|y|(2b+1)}{2}+\frac{|y|\sqrt{|D|}}{2}.
\end{align}
Taking the minimum over all primitive elements $\alpha \in I$, we see that
\begin{equation}
\label{N3-2}
\min \{ h(\alpha) : \alpha \in I, K = \que(\alpha) \} \leq \frac{(2b+1) + \sqrt{|D|}}{2},
\end{equation}
where the inequality is obtained by taking $x=0$, $y=1$. Putting together \eqref{N_h}, \eqref{N1-2}, \eqref{N2-2} and \eqref{N3-2}, we obtain the $D \equiv 1 (\md 4)$ case of Theorem~\ref{quad_bound} in case $g=1$.
\medskip

\proof[Proof of Theorem~\ref{quad_bound}]
Let $I$ be an ideal with the canonical basis $\{a,b+g\delta\}$ and $J = \frac{1}{g} I$. Then for any $\alpha \in J$ and the corresponding $g\alpha \in I$,
\begin{equation}
\label{hgI}
h(g\alpha) = \left( \prod_{j=1}^2 \max \{1,|\sigma_i(g \alpha)|\} \right)^{\frac{1}{2}} \leq g \left( \prod_{j=1}^2 \max \{1,|\sigma_i(\alpha)|\} \right)^{\frac{1}{2}} = g h(\alpha).
\end{equation}
Further, $\Nn_K(I) = g \Nn_K(J)$. Take $\alpha \in J$ be a primitive element of bounded height as obtained above in Cases 1 and 2, then $g\alpha \in I$ is also a primitive element and the result follows.
\endproof

\bigskip

\section{Primitive elements with avoidance conditions}
\label{pf_ideal_miss}

The goal of this section is to prove Theorem~\ref{ideal_miss}. Our main tool is Theorem~\ref{dioph_avoid}, so we will set things up to apply it. First recall that a union of ideals $\bigcup_{i=1}^s J_k$ is an ideal if and only if there exists some $1 \leq i \leq s$ such that $J_1,\dots,J_s \subseteq J_i$: if this is not the case, the union would not be closed under addition. Since they are all properly contained in $I$, we conclude that $I \neq \bigcup_{i=1}^s J_i$, and so there exists $\alpha \in I \setminus \bigcup_{i=1}^s J_i$. Let us then define lattices
$$\Omega = \Sigma_K(I),\ \Lambda_i = \Sigma_K(J_i)\ \forall\ 1 \leq i \leq s,\ \Lambda = \Sigma_K(J)$$
in the Euclidean space $K_{\real}$. Let
\begin{equation}
\label{P_poly}
P(x_1,\dots,x_d) = \prod_{1 \leq i < j \leq d} (x_i - x_j) \in \zed[x_1,\dots,x_d],
\end{equation}
and notice that an element $\alpha \in K$ is primitive if and only if $P(\Sigma_K(\alpha)) \neq 0$. Indeed, an element $\alpha \in K$ is primitive if and only if all of its algebraic conjugates $\sigma_1(\alpha),\dots,\sigma_d(\alpha)$ are distinct, but these are precisely the coordinates of the vector $\Sigma_K(\alpha)$ in $K_{\real}$. Degree of this polynomial $P$ is $m = \binom{d}{2}$, and so Theorem~\ref{dioph_avoid} guarantees the existence of a point $\bz \in \Omega \setminus \left( \bigcup_{i=1}^s \Lambda_i \right)$
such that $P(\bz) \neq 0$ and~\eqref{z-bound} is satisfied. Let $\alpha \in I$ be such that $\bz = \Sigma_K(\alpha)$. We now want to rewrite the inequality~\eqref{z-bound} in terms of the invariants of the ideals and the number field $K$ and use it to estimate $h(\alpha)$.
\smallskip

Since $D_i$ and $D$ are indices of $\Lambda_i$ and $\Lambda$, respectively, in $\Omega$, we have
$$D_i = \frac{\det \Lambda_i}{\det \Omega} = \frac{\Nn_K(J_i)}{\Nn_K(I)},\ D = \frac{\det \Lambda}{\det \Omega} = \frac{\Nn_K(J)}{\Nn_K(I)},$$
by~\eqref{ideal_det}, and $\Delta = \det \Omega = \Nn_K(I) |\Delta_K|^{1/2}$. Now, the ``cube"
$$C_d(1) = \left\{ \bx \in K_{\real} : |\bx| \leq 1 \right\}$$
is the Cartesian product of $r_1$ intervals $[-1,1]$ and $r_2$ circles of radius $1$. Hence, $C_d(1)$ is a convex $\bo$-symmetric set with $d$-dimensional volume
$$\Vol_d(C_d(1)) = 2^{r_1} \pi^{r_2}.$$
Finally, notice that
\begin{eqnarray*}
1 & = & \min \left\{  \left( \prod_{v \mid \infty} |\beta|_v \right)^{1/d} : \beta \in \O_K \right\} \\
& \leq & \min \left\{ \max_{1 \leq j \leq d} |\sigma_j(\beta)| : \beta \in J \setminus \{0\} \right\} = \lambda_1(\Lambda) \leq \Nn_K(J),
\end{eqnarray*}
since $\Nn_K(J) \in J \cap \zed_{>0}$, and so $|\sigma_j(\Nn_K(J))| = \Nn_K(J)$ for every $1 \leq j \leq d$. Putting all of the above observations together with Lemma~\ref{ht_sigma}, we obtain~\eqref{main}. This completes the proof of Theorem~\ref{ideal_miss}.

\bigskip

\section{Positive points with avoidance conditions}
\label{cone}

In this section, we let our Euclidean space $\EE^d$ to be just $\real^d$. We consider ``small" {\it positive} points in a lattice outside of a union of sublattices, at which a given polynomial does not vanish. A previous result in~\cite{lf_siki} provided a bound on the smallest norm of a nonzero positive point in a lattice, i.e., a point in the intersection of a lattice with the positive orthant in~$\real^d$. Here, we extend this result to include additional avoidance conditions. Specifically, we prove the following theorem.

\begin{thm} \label{pos_cone} Let $\Omega \subseteq \real^d$ be a lattice of full rank and determinant $\Delta$, $\Lambda_1,\dots,\Lambda_s \subseteq \Omega$ its sublattices with indices $[\Omega:\Lambda_i] = D_i$ for each $1 \leq i \leq s$, and $\Lambda = \bigcap_{i=1}^s \Lambda_i$ be a sublattice of $\Omega$ of index $D \leq D_1 \cdots D_s$. Assume that $\Omega \not\subseteq \bigcup_{i=1}^s \Lambda_i$. Let $P(x_1,\dots,x_d) \in \real[x_1,\dots,x_d]$ be a polynomial of degree $m$. Then there exists a point $\bz \in \Omega \setminus \bigcup_{i=1}^s \Lambda_i$ so that
$$P(\bz) \neq 0,\ z_i \geq 0\ \forall\ 1 \leq i \leq d,$$
and
$$|\bz| < D (m + 2 ) (\mu(\Lambda)+1) \left( \frac{D \Delta}{\lambda_1(\Lambda)^{d-1}} \left( \sum_{i=1}^s \frac{1}{D_i} - \frac{s - 1}{D} + \frac{\lambda_1(\Lambda)^d}{D\Delta} \right) + \sum_{i=1}^d \lambda_i(\Lambda) \right),$$
where $\lambda_i(\Lambda)$ are the successive minima of $\Lambda$ with respect to the cube $C_d(1)$ as defined in~\eqref{sm} and
$$\mu(\Lambda) = \min \left\{ T \in \real_{> 0} : B_d(T) + \Lambda = \real^d \right\}$$
is the covering radius of $\Lambda$ with respect to the unit ball $B_d(1)$. 
\end{thm}

\proof
Let us write
$$\Omega^+ = \left\{ \bx \in \Omega : x_j \geq 0\ \forall\ 1 \leq j \leq d \right\},$$
then $\Lambda^+ \subset \Omega^+$. The restricted successive minima of $\Lambda^+$, defined as
$$\lambda_i(\Lambda^+) := \min \left\{ T \in \real_{> 0} : \dim_{\real} \spn_{\real} \left( \Lambda^+ \cap C_d(T) \right) \geq i \right\},$$
were studied in~\cite{lf_siki}. Theorem~1.2 of~\cite{lf_siki} established that
\begin{equation}
\label{smL+}
\lambda_1(\Lambda^+) \leq 2\mu(\Lambda)+1,\ \lambda_i(\Lambda^+) \leq 2\lambda_i(\Lambda) (\mu(\Lambda) + 1)\ \forall\ 2 \leq i \leq d.
\end{equation}
By Lemmas~3.1 and~3.2 of~\cite{lf_siki}, there exist linearly independent vectors $\bv_1,\dots,\bv_d \in \Lambda^+$ such that
$$|\bv_i| \leq \lambda_i(\Lambda^+)\ \forall\ 1 \leq i \leq d, \text{ and } v_{1j} \geq 1\ \forall\ 1 \leq j \leq d.$$
Define the vectors $\bu_1 = \bv_1$ and $\bu_i = \bv_i + \bv_1$ for every $2 \leq i \leq d$, then
$$\bu_1,\dots,\bu_d \in \Lambda^+ \text{ and } u_{ij} \geq 1\ \forall\ 1 \leq i,j, \leq d.$$
Further, $|\bu_1| \leq 2\mu(\Lambda)+1$ and for every $2 \leq i \leq d$,
\begin{equation}
\label{u_vectors}
|\bu_i| \leq |\bv_i| + |\bv_1| \leq 2 (\lambda_i(\Lambda) + 1) (\mu(\Lambda) + 1) - 1.
\end{equation}
Let $\bx \in \Omega \setminus \bigcup_{i=1}^s \Lambda_i$ be as in~\eqref{ht}. As in our argument in Section~\ref{DA}, there must exist a subset of $d-1$ vectors among $\bu_1,\dots,\bu_d$ with which $\bx$ is linearly independent, assume these are $\bu_2,\dots,\bu_d$, since they have larger norm. Notice that the vector
$$\bwy := \bx + |\bx| \bu_1 \in \Omega^+,$$
since its coordinates are of the form
$$y_i = x_i + |\bx| v_{1i} \geq 0.$$
On the other hand, $\bwy \not\in \bigcup_{i=1}^s \Lambda_i$, since $\bx$ is not contained in any $\Lambda_i$ while $\bu_1$ is contained in each of them. Additionally,
\begin{equation}
\label{y_bnd}
|\bwy| \leq |\bx| (1+|\bu_1|) \leq 2|\bx| (\mu(\Lambda)+1).
\end{equation}
Let $P(x_1,\dots,x_d) \in \real[x_1,\dots,x_d]$ be a polynomial of degree $m \geq 1$. We now argue as in Section~\ref{DA}. Let the sets $S_1$ and $S_2$ be as in~\eqref{S1S2}, then~\eqref{p_avoid} implies that there exists a vector 
$$\bxi = (\xi_1,\xi_2,\dots,\xi_d) \in S_1 \times S_2 \times \dots \times S_2$$
so that 
$$P\left( \xi_1 \bwy + \sum_{i=2}^d \xi_i \bu_i \right) \neq 0.$$
Let $\bz = \xi_1 \bwy + \sum_{i=2}^d \xi_i \bu_i \in \Omega^+ \setminus \bigcup_{i=1}^s \Lambda_i$ for this choice of $\bxi$, then by~\eqref{u_vectors} and~\eqref{y_bnd}, we have
$$|\bz| \leq |\bxi| \left( |\bwy| + \sum_{i=2}^d |\bu_i| \right) \leq D (m + 2 ) (\mu(\Lambda)+1) \left( |\bx| + \sum_{i=1}^d \lambda_i(\Lambda) \right).$$
Combining this last observation with~\eqref{ht} and taking into account that $\Vol_d(C_d(1)) = 2^d$ finishes the proof of the theorem.
\endproof
\smallskip

We can now apply this theorem to the number field situation. Let $K$ be a totally real number field of degree $d$ and discriminant~$\Delta_K$ with embeddings
$$\sigma_1,\dots,\sigma_d : K \hookrightarrow \real,$$
which define the Minkowski embedding $\Sigma_K = (\sigma_1,\dots,\sigma_d) : K \hookrightarrow \real^d$. Let $I \subseteq \O_K$ be an ideal and write $I^+$ for the additive semigroup of totally positive elements in $I$, i.e.
$$I^+ = \left\{ \alpha \in I : \sigma_i(\alpha) \geq 0\ \forall\ 1 \leq i \leq d \right\}.$$
Same as in Section~\ref{pf_ideal_miss}, let $J_1,\dots,J_s$ be distinct ideals properly contained in $I$ and let $J = J_1 \cdots J_s$. Let
$$\Omega = \Sigma_K(I),\ \Lambda_i = \Sigma_i(J_i)\ \forall\ 1 \leq i \leq s,\ \Lambda = \Sigma_K(J)$$
be lattices in $\real^d$. Let
$$P(x_1,\dots,x_d) = \prod_{1 \leq i < j \leq d} (x_i - x_j) \in \zed[x_1,\dots,x_d],$$
so an element $\alpha \in K$ is primitive if and only if $P(\Sigma_K(\alpha)) \neq 0$. Notice that $\Sigma_K(I^+) = \Sigma_K(I)^+ = \Omega^+$. Now, Lemmas~4.1 and~4.2 of~\cite{lf_siki} combined guarantee the existence of $\que$-linearly independent elements $a_1,\dots,a_d \in J$ such that
\begin{equation}
\label{I_sm-1}
\sum_{i=1}^d \lambda_i(\Lambda) \leq \sum_{i=1}^d h(a_i)^d \leq d \prod_{i=1}^d h(a_i)^d  \leq d \left( \Nn_K(J) \sqrt{ |\Delta_K|} \right)^d.
\end{equation}
Further, Lemmas~4.2 of~\cite{lf_siki} asserts that the covering radius of $\Lambda$ satisfies the inequality
\begin{equation}
\label{I_sm-2}
\mu(\Lambda) \leq \frac{d^{3/2}}{2} \Nn_K(J) \sqrt{|\Delta_K|}.
\end{equation}
Putting these observations together with Theorem~\ref{pos_cone} and expressing indices and determinants as in  Section~\ref{pf_ideal_miss}, we obtain the following corollary.

\begin{cor} \label{ideal_pos} There exists a primitive totally positive element $\alpha \in I \setminus \bigcup_{i=1}^s J_i$ so that
\begin{eqnarray*}
h(\alpha) & \leq & \left( \binom{d}{2} + 2 \right)  \left( \frac{d^{3/2}}{2} \Nn_K(J) \sqrt{|\Delta_K|} +1 \right) \frac{\Nn_K(J)^2 \sqrt{|\Delta_K|}}{\Nn_K(I)} \times \\
& \times & \left( \sum_{i=1}^s \frac{\Nn_K(I)}{\Nn_K(J_i)} - \frac{(s - 1) \Nn_K(I)}{\Nn_K(J)} + \frac{\Nn_K(J)^{d-1}}{\sqrt{|\Delta_K|}} + d \left( \Nn_K(J) \sqrt{ |\Delta_K|} \right)^{d-1} \right).
\end{eqnarray*}
\end{cor}

\bigskip

\section{Non-sparse polynomials with bounded Mahler measure}
\label{mahler}

Let $K$ be a number field of degree $d$, then~\eqref{wd1} guarantees that there exists a primitive element $\alpha \in \O_K$ with $h(\alpha) \leq |\Delta_K|^{\frac{1}{d}}$. Let 
$$f_{\alpha}(x) = x^d + \sum_{k=0}^{d-1} a_k x^k \in \zed[x]$$
be the minimal polynomial of $\alpha$. Then $f_{\alpha}(x)$ is a monic irreducible polynomial with integer coefficients and $K \cong \que[x]/\left< f_{\alpha}(x) \right>$. Let $\sigma_1,\dots,\sigma_d$ be the embeddings of $K$, as usual, ordered in such a way that $\sigma_1(\alpha) = \alpha$. Then 
$$\alpha_k = \sigma_k(\alpha_1),\ \forall\ 1 \leq k \leq d$$
are all the roots of $f_{\alpha}(x)$ with $\alpha_1 = \alpha$. Since $f_{\alpha}(x)$ is monic, all of these roots are algebraic integers. The Mahler measure of $f_{\alpha}(x)$ is given by
$$\M(f_{\alpha}) = \prod_{k=1}^d \max \{1, |\alpha_k| \} = h(\alpha)^d \leq |\Delta_K|.$$
Hence, the result of Pazuki and Widmer~\cite[Lemma~7.1]{widmer} (see equation~\eqref{wd0} above) can be reformulated to say that there exists a monic irreducible polynomial $f(x) \in \zed[x]$ such that $K \cong \que[x]/\left< f(x) \right>$ and $\M(f) \leq |\Delta_K|$. We can use our Diophantine avoidance method to obtain a similar result with additional non-vanishing conditions.

\begin{thm} \label{mm} Given a number field $K$ of degree $d$, there exists a monic irreducible polynomial $f(x) \in \zed[x]$ with all nonzero coefficients such that $K \cong \que[x]/\left< f(x) \right>$ and 
$$\M(f) \leq \left\{ \left( \frac{4}{\pi} \right)^{r_2}  \left( \frac{d(d^2-d+2)}{2} \right)  |\Delta_K|^{1/2} \right\}^d.$$
\end{thm}

\proof
Our argument here is similar to the proof of Theorem~\ref{ideal_miss}. Let the lattice $\Omega = \Sigma_K(\O_K)$ in $\EE^d = K_{\real}$ be the image of $\O_K$ under the Minkowski embedding. For $1 \leq k \leq d$, let $e_k(x_1,\dots,x_d)$ be the elementary symmetric polynomial of degree $k$. Let $\alpha \in \O_K$ and $\baa := \Sigma_K(\alpha) = (\alpha_1,\dots,\alpha_d) \in \Omega$. Then $\alpha$ is a primitive element in $K$ if and only if its minimal polynomial is of the form
$$f_{\alpha}(x) = x^d + \sum_{k=1}^{d} e_k(\alpha_1,\dots,\alpha_d) x^{d-k},$$
which is equivalent to the condition that $K \cong \que[x]/\left< f_{\alpha}(x) \right>$. Let
$$Q(x_1,\dots,x_d) = P(x_1,\dots,x_d) \prod_{k=1}^{d-1} e_k(x_1,\dots,x_d),$$
where $P(x_1,\dots,x_d)$ is as in~\eqref{P_poly}, then
$$\deg Q = \deg P + \sum_{k=1}^{d-1} k = d(d-1).$$
Notice then that $\alpha$ is primitive and $f_{\alpha}(x)$ has all nonzero coefficients if and only if $Q(\baa) \neq 0$. Hence, we want to find $\baa \in \Omega$ such that $Q(\baa) \neq 0$ and $|\baa|$ bounded. Let
$$\bv_1,\dots,\bv_d \in \Omega$$
be vectors corresponding to the successive minima of $\Omega$ with respect to the cube $C_d(1) \subset K_{\real}$. Then $\Vol_d(C_d(1)) = 2^{r_1} \pi^{r_2}$ and $\det \Omega = |\Delta_K|^{1/2}$, and so Minkowski's successive minima theorem implies that
\begin{equation}
\label{mink-1}
\lambda_d(\Omega) \leq \frac{2^d |\Delta_K|^{1/2}}{\lambda_1(\Lambda)^{d-1} \Vol_d(C_d(1))} \leq \left( \frac{4}{\pi} \right)^{r_2} |\Delta_K|^{1/2},
\end{equation}
since 
$$\lambda_1(\Omega) = \min \left\{ \max_{1 \leq k \leq d} |\sigma_k(\alpha)| : \alpha \in \O_K \right\} = 1.$$
Let $S = \left\{ - \left[ \frac{d(d-1)}{2} \right] - 1, \dots,\left[ \frac{d(d-1)}{2} \right] + 1 \right\}$, and for each vector $\bxi \in S^d$ define
$$\baa(\bxi) = \sum_{k=1}^d \xi_k \bv_k.$$
Then Theorem~4.2 of~\cite{lf3} implies that there exists $\bxi \in S^d$ such that
$$P(\baa(\bxi)) \neq 0.$$
Let $\alpha \in \O_K$ be such that $\baa(\bxi) = \Sigma_K(\alpha)$ for this choice of $\bxi$. Then Lemma~\ref{ht_sigma} together with~\eqref{mink-1} implies that
$$h(\alpha) \leq |\baa(\bxi)| \leq d \left( \frac{d(d-1)}{2} + 1 \right) \lambda_d(\Omega) \leq \left( \frac{4}{\pi} \right)^{r_2}  \left( \frac{d(d^2-d+2)}{2} \right)  |\Delta_K|^{1/2}.$$
If $f(x)$ is the minimal polynomial of this $\alpha$, then $\M(f) \leq h(\alpha)^d$, and the result follows.
\endproof

It is interesting to remark that some previous literature featured lower bounds on the Mahler measure of {\it fewnomials}, i.e., polynomials with few nonzero coefficients; see, for instance,~\cite{dobrowolski} and references within, such as a classical paper of Mahler~\cite{mahler}. Our Theorem~\ref{mm} goes in the opposite direction, providing an upper bound on the Mahler measure of polynomials with all nonzero coefficients generating a given number field. For comparison, inequality (7) of~\cite{mahler} gives an upper bound on $\M(f)$, however it is in terms of height (maximum of absolute values of coefficients) of the polynomial~$f$, whereas we prove the existence of a non-sparse generating polynomial~$f$ for~$K$ with Mahler measure bounded in terms of the invariants of $K$.
\bigskip

\section{Small-height ideal generators}
\label{small_gen}

In this section, we prove Theorem~\ref{ideal_gen}. As in Section~\ref{quad}, we consider a quadratic number field $K = \que(\sqrt{D})$ for squarefree integer $D$, and let $I \subseteq \O_K$ be an ideal with the canonical basis $a,b+g\delta$. Suppose $I = \left< \mu \right>$ is a principal ideal, then $\Nn_K(\mu) = \Nn_K(I)$. Also recall that we are defining the quantity $H(I)$, the ``height" of the ideal $I$ as in~\eqref{HI}.
\smallskip

{\it Case 1.} Assume $D \not\equiv 1\ (\md 4)$, then $\Nn_K(\mu) = \Nn_K(I) = ag$. On the other hand,
$$\mu = ax+(b-g\sqrt{D})y = (ax+by)-yg\sqrt{D},$$
for some integers $x,y$ such that $y \neq 0$, and thus
$$ag = \Nn_K(\mu) = (ax+by)^2-Dgy^2 = a^2x^2+2abxy+(b^2-Dg)y^2.$$
This implies that the quadratic equation with integer coefficients
\begin{equation}
\label{f_I}
f_I(x,y) := \frac{a}{g} x^2+\frac{2b}{g} xy+\left( \frac{b^2-Dg}{ag} \right) y^2 = 1
\end{equation}
has a solution in integers $x,y$. Define
$$|f_I| := \max \left\{1, \frac{|a|}{g},\frac{|2b|}{g},\frac{|b^2-D|}{ag} \right\} = H(I)$$
to be the maximum of absolute values of the coefficients of $f_I$.  Then a theorem of Kornhauser~\cite[Theorem~1]{kornhauser} guarantees that~\eqref{f_I} has an integer solution $(x,y)$ with
\begin{equation}
\label{kornh}
\max \{|x|,|y|\} \leq (14|f_I|)^{5|f_I|}.
\end{equation}
Now, using inequalities analogous to \eqref{h_a_mu} and applying~\eqref{kornh} yields the inequality
\begin{equation}
\label{hmu-1}
h(\mu) \leq |x|a + |y|b + |y| g \sqrt{|D|} \leq \left( a+b+g\sqrt{|D|} \right) (14|f_I|)^{5|f_I|}.
\end{equation}
\smallskip

{\it Case 2.} Assume $D \equiv 1\ (\md 4)$, then $\Nn_K(\mu) = \Nn_K(I) = \frac{ag}{2}$. On the other hand,
$$\mu = ax+\left(b+g \left( \frac{1-\sqrt{D}}{2} \right) \right)y = \left(ax+by+\frac{gy}{2}\right)-y\left(\frac{g\sqrt{D}}{2}\right)$$
for some integers $x,y$ such that $y \neq 0$, and thus
\begin{eqnarray*}
\frac{ag}{2} & = & \Nn_K(\mu) = \left( ax+by+\frac{yg}{2} \right)^2-\frac{gD}{4}y^2 \\
& = & a^2x^2+a(2b+g)xy+\left( b^2+bg-\frac{D-1}{4} g \right)y^2.
\end{eqnarray*}
This implies that the quadratic equation with integer coefficients
\begin{equation}
\label{f_I-1}
f_I(x,y) := \frac{2a}{g} x^2 + \frac{2(2b+g)}{g} xy+\left( \frac{4b^2+4b-Dg+g}{2ag} \right) y^2 = 1
\end{equation}
has a solution in integers $x,y$. Define
$$|f_I| := \max \left\{1, \frac{|2a|}{g}, \frac{2 |2b+g|}{g} ,\frac{|4b^2+4b-Dg+g|}{2ag} \right\} = H(I)$$
to be the maximum of absolute values of the coefficients of $f_I$. Now, using inequalities analogous to \eqref{h_a_mu-1} and applying~\eqref{kornh} yields the inequality
\begin{equation}
\label{hmu-2}
h(\mu) \leq |x|a+\frac{|y|(2b+g)}{2}+\frac{|y|g\sqrt{|D|}}{2} \leq \left( a + \frac{2b+g+ g\sqrt{|D|}}{2} \right) (14|f_I|)^{5|f_I|}.
\end{equation}
The theorem now follows upon combining~\eqref{hmu-1} and~\eqref{hmu-2}. Notice that the element $\mu$ we produced here is primitive since $y \neq 0$ in both cases.
\bigskip

\noindent
{\bf Acknowledgement:} We wish to thank the anonymous referee for a thorough reading and multiple suggestions that improved the quality of this paper.

\bigskip

%\nocite{*}
\bibliographystyle{plain}  % Here the bibliography 

\end{document}